# A New Class of Gamma distribution

Cícero Carlos Ramos de Brito

*Statistics and Informatics Department, Universidade Federal Rural de Pernambuco, Recife, Brazil.*

Address: Av. Professor Cláudio Selva, 178, Dois Irmãos, Recife, Pernambuco, Brazil.

CEP: 52171-260

cicerocarlosbrito@yahoo.com.br

Leandro Chaves Rêgo

*Department of Statistics, Universidade Federal de Pernambuco, Recife, Brazil*

leandro@de.ufpe.br

Wilson Rosa de Oliveira

*Statistics and Informatics Department, Universidade Federal Rural de Pernambuco, Recife, Brazil.*

wilson.rosa@gmail.com

# A New Class of Gamma distribution

This paper presents a new class of probability distributions generated from the gamma distribution. For the new class proposed, we present several statistical properties, such as the risk function, the density expansions, Moment-generating function, characteristic function, the moments of order m, central moments of order m, the log likelihood and its partial derivatives and also entropy, kurtosis, symmetry and variance. These same properties are determined for a particular distribution within this new class that is used to illustrate the capability of the proposed new class through an application to a real data set. The database presented in Choulakian and Stephens (2001) was used. Six models are compared and for the selection of these models were used the Akaike Information Criterion (AIC), the Akaike Information Criterion corrected (AICc), Bayesian Information Criterion (BIC), Hannan Quinn Information Criterion (HQIC) and tests of Cramer-Von Mises and Anderson-Darling to assess the models fit. Finally, we present the conclusions from the analysis and comparison of the results obtained and the directions for future work.

## 1. Introduction

The Gamma distribution is used in a variety of applications including queue, financial and weather models. It can naturally be considered as the distribution of the waiting time between events distributed according to a Poisson process. It is a biparamétrica distribution, whose density is given by:

$$f(t) = \frac{\beta^\alpha}{\Gamma(\alpha)} t^{\alpha-1} e^{-\beta t}, t > 0,$$

where $\alpha > 0$ is a shape parameter and $\beta > 0$ is the reciprocal of a scale parameter.

Due to the importance of this distribution recently some new distributions as well as families of probability distributions based on generalizations of the Gamma distribution have been proposed. The first is based on the family of exponentiated distribution defined by Mudholkar *et al.* (1995). Given a distribution with continuous distribution function $F$ d its *generalization* or *exponentialization* $G(x)$ is obtained by $F(x) = G^a(x)$, with $a > 0$. Gupta *et al.* (1998) proposed and studied some properties exponentiated Gamma distribution.

Cordeiro *et al.* (2011) extended the exponentiated Gamma distribution defining a new distribution called Exponentiated Generalized Gamma Distribution with four parameters, which is capable of modeling bathtub shaped failure rate phenomena.

Zografos and Balakrishnan (2009) defined a family of probability distributions based on the integration of a Gamma distribution as follows:

$$F(x) = \frac{\beta^\alpha}{\Gamma(\alpha)} \int_0^{-\ln(1-G(x))} t^{\alpha-1} e^{-\beta t} dt,$$

where $G(x)$ is an arbitrary distribution function. When $\alpha = n + 1$ and $\beta = 1$ this distribution coincides with the distribution of the *n-th* highest value record (ALZAATREH *et al.* 2014).

A new family of probability distributions which is also based on the integration of the Gamma distribution has been proposed by Silva (2013). This author defined this new family as follows:

$$F(x) = 1 - \frac{\beta^\alpha}{\Gamma(\alpha)} \int_0^{-\ln(G(x))} t^{\alpha-1} e^{-\beta t} dt,$$

where $G(x)$ is an arbitrary distribution function. When $\alpha = n + 1$ and $\beta = 1$ this distribution coincides with the distribution of the *n*th smallest value record (ALZAATREH *et al.* 2014).

Following the line of work of Zografos and Balakrishnan (2009) and Silva (2013), our goal in this work is to propose a new family of distributions based on Gamma distribution. The family of distributions proposed here is the following:

$$H_{G_1}(x) = \int_{\frac{1-G_1(x)}{G_1(x)}}^{+\infty} \frac{\beta^\alpha}{\Gamma(\alpha)} t^{\alpha-1} e^{-\beta t} dt,$$

where $G_1(x)$ is an arbitrary distribution function and $H_{G_1}(x)$ and has the same support as the distribution $G_1(x)$. We shall call this new class $\frac{1-G_1}{G_1}$ Gama Class. The statistical properties of this new class, such as mean, variance, standard deviation, mean deviation, kurtosis, skewness, moment generating function, characteristic function and graphical analysis, are derived.

Then, to illustrate the applicability of the proposed new family, we consider the particular case of the distribution obtained when considering that $G_1(x)$ is the distribution function of an exponential random variable. The statistical properties of this new distribution are also derived and to illustrate its potentiality, an application to a set of real data is performed. For this, we used the database presented in the work Choulakian and Stephens (2001) to see if the models are well adjusted to this data. As comparative criteria of fitness of the models, it was considered: the Akaike (AIC) (AKAIKE, 1972), the Akaike Fixed (AIC) (BURNHAM AND ANDERSON, 2002), the Bayesian information criterion (BIC) (SCHWARTZ, 1978), the Hannan-Quinn information criterion (HQIC) (HANNAN AND QUINN, 1979), and the Cramer-von Mises (DARLING, 1957) and Anderson-Darling (ANDERSON and DARLING, 1952) tests. Both hypothesis tests, Anderson-Darling and Cramér-von Mises, are discussed in detail by Chen and Balakrishnan (1995) and belong to the class of quadratic statistics based on the empirical distribution function, because they work with the squared differences between the empirical distribution and the hypothetical.

This paper is organized as follows. In Section 2, we describe the statistical properties of the proposed Gamma class. In section 3, we studied a special case of this new class in the case in which the distribution function is an exponential random variable. We study the properties of this new distribution and perform an application to real data illustrating its potential. Finally, we present our conclusions in Section 4.

**2. Proposed Model**

*2.1 Obtaining a class of probability distributions*

*2.1.1 Model functional class range (1-$G_1$) / $G_1$*

The $\frac{1-G_1}{G_1}$ Gamma class is defined by the cumulative distribution function:

$$H_{G_1}(x) = \int_{\frac{1-G_1(x)}{G_1(x)}}^{+\infty} \frac{\beta^\alpha}{\Gamma(\alpha)} t^{\alpha-1} e^{-\beta t} dt,$$

which is equivalent to

$$H_{G_1}(x) = 1 - \int_0^{\frac{1-G_1(x)}{G_1(x)}} \frac{\beta^\alpha}{\Gamma(\alpha)} t^{\alpha-1} e^{-\beta t} dt.$$

If the distribution $G_1(x)$ has density $g_1(x)$ the class will have a probability density function given by

$$h_{G_1}(x) = \frac{g_1(x)}{G_1^2(x)} \frac{\beta^\alpha}{\Gamma(\alpha)} \left(\frac{1-G_1(x)}{G_1(x)}\right)^{\alpha-1} \exp\left(-\beta\left(\frac{1-G_1(x)}{G_1(x)}\right)\right).$$

Such functions can be rewritten as a sum of *exponentiated* distributions, as follows. As

$$\exp\left(-\beta \frac{1-G_1(x)}{G_1(x)}\right) = \sum_{k=0}^{\infty} \frac{(-1)^k \beta^k}{k!} G_1^{-k}(1-G_1(x))^k,$$

we have to

$$h_{G_1}(x) = (g_1(x)) \frac{\beta^\alpha}{\Gamma(\alpha)} \sum_{k=0}^{\infty} \frac{(-1)^k \beta^k}{k!} G_1^{-\alpha-k-1}(1-G_1(x))^{\alpha+k-1}.$$

Furthermore, as

$$(1-G_1(x))^{k+\alpha-1} = \sum_{j=0}^{\infty} \binom{k+\alpha-1}{j} (-1)^j G_1^j(x),$$

it follows that

$$h_{G_1}(x) = \sum_{k=0}^{\infty} \sum_{j=0}^{\infty} \frac{(-1)^{k+j} \beta^{\alpha+k}}{k! \Gamma(\alpha)} \binom{k+\alpha-1}{j} g_1(x) G_1^{j-\alpha-k-1}(x)$$

Since $H_{G_1}(x) = \int_{-\infty}^{x} h_{G_1}(t)dt$, we can rewrite the distribution function as

$$H_{G_1}(x) = \sum_{k=0}^{\infty} \sum_{j=0}^{\infty} \frac{(-1)^{k+j} \beta^{\alpha+k}}{k! \Gamma(\alpha)} \binom{k+\alpha-1}{j} \int_{-\infty}^{x} g_1(t) G_1^{j-\alpha-k-1}(t)dt.$$

Therefore,

$$H_{G_1}(x) = \sum_{k=0}^{\infty} \sum_{j=0}^{\infty} \frac{(-1)^{k+j} \beta^{\alpha+k}}{k!(j-\alpha-k)\Gamma(\alpha)} \binom{k+\alpha-1}{j} G_1^{j-\alpha-k}(x).$$

If the distribution $G_1(x)$ is discrete, $H_{G_1}(x)$ is also discrete and we have that $P(X = x_l) = F(x_l) - F(x_{l-1})$. Therefore,

$$P(X = x_l) = \sum_{k=0}^{\infty} \sum_{j=0}^{\infty} \frac{(-1)^{k+j} \beta^{\alpha+k}}{k!(j-\alpha-k)\Gamma(\alpha)} \binom{k+\alpha-1}{j} \left(G_1^{j-\alpha-k}(x_l) - G_1^{j-\alpha-k}(x_{l-1})\right). \blacksquare$$

In addition, we can obtain the risk function of the new $\frac{1-G_1}{G_1}$Gamma class as follows:

$$R_{G_1}(x) = \frac{h_{G_1}(x)}{1 - H_{G_1}(x)}$$

$$R_{G_1}(x) = \frac{\frac{g_1(x)}{G_1^2(x)} \frac{\beta^\alpha}{\Gamma(\alpha)} \left(\frac{1-G_1(x)}{G_1(x)}\right)^{\alpha-1} \exp\left(-\beta\left(\frac{1-G_1(x)}{G_1(x)}\right)\right)}{\int_0^{\frac{1-G_1(x)}{G_1(x)}} \frac{\beta^\alpha}{\Gamma(\alpha)} t^{\alpha-1} e^{-\beta t} dt}.$$

∎

Using the density and distribution function expansions, we can get the statistical properties of the new class, as discussed below.

*2.1.2 Expansion for the moments of order m for the (1-$G_1$)/$G_1$ Gamma Class*

The following is the development of the expansion calculations for the moments of order $m$ for the $\frac{1-G_1}{G_1}$ Gamma Class. As

$$\mu_m = E(X^m) = \int_{-\infty}^{+\infty} x^m dF(x),$$

we have to

$$\mu_m = \int_{-\infty}^{+\infty} x^m \sum_{k=0}^{\infty} \sum_{j=0}^{\infty} \frac{(-1)^{k+j} \beta^{\alpha+k}}{k!\Gamma(\alpha)} \binom{k+\alpha-1}{j} g_1(x) G_1^{j-\alpha-k-1}(x) dx.$$

Therefore,

$$\mu_m = \sum_{k=0}^{\infty} \sum_{j=0}^{\infty} \frac{(-1)^{k+j} \beta^{\alpha+k}}{k!\Gamma(\alpha)} \binom{k+\alpha-1}{j} \tau_{m,0,j-\alpha-k-1},$$

*where*

$$\tau_{m,\eta,r} = E(X^m g_1(X)^\eta G_1(X)^r) = \int_{-\infty}^{+\infty} x^m g_1(x)^\eta G_1(x)^r dG_1(x).$$

∎

In particular, we have the following expansion of the mean for the $\frac{1-G_1}{G_1}$ Gamma Class

$$\mu = \mu_1 = \sum_{k=0}^{\infty} \sum_{j=0}^{\infty} \frac{(-1)^{k+j} \beta^{\alpha+k}}{k!\Gamma(\alpha)} \binom{k+\alpha-1}{j} \tau_{1,0,j-\alpha-k-1}.$$

∎

*2.1.3 Expansion for the moment generating function for the (1-$G_1$)/$G_1$ Gamma Class.*

The following is the development of the expansion calculations for the moment generating function for the $\frac{1-G_1}{G_1}$ Gamma Class. As

$$M_X(t) = E(e^{tX}) = \int_{-\infty}^{+\infty} e^{tx} dF(x),$$

we have

$$M_X(t) = \sum_{k=0}^{\infty} \sum_{j=0}^{\infty} \frac{(-1)^{k+j} \beta^{\alpha+k}}{k!\,\Gamma(\alpha)} \binom{k+\alpha-1}{j} \int_{-\infty}^{+\infty} e^{tx} g_1(x) G_1^{j-\alpha-k-1}(x) dx.$$

Using the fact that

$$e^{tx} = \sum_{m=0}^{\infty} \frac{t^m x^m}{m!},$$

we can rewrite

$$M_X(t) = \sum_{k=0}^{\infty} \sum_{j=0}^{\infty} \sum_{m=0}^{\infty} \frac{(-1)^{k+j} \beta^{\alpha+k} t^m}{k!\,m!\,\Gamma(\alpha)} \binom{k+\alpha-1}{j} \int_{-\infty}^{+\infty} x^m g_1(x) G_1^{j-\alpha-k-1}(x) dx.$$

Therefore,

$$M_X(t) = \sum_{k=0}^{\infty} \sum_{j=0}^{\infty} \sum_{m=0}^{\infty} \frac{(-1)^{k+j} \beta^{\alpha+k} t^m}{k!\,m!\,\Gamma(\alpha)} \binom{k+\alpha-1}{j} \tau_{m,0,j-\alpha-k-1}$$

∎

Similarly, one can establish the following expansion for the characteristic function for the (1-G1)/G1 Gamma Class.

$$\varphi_X(t) = \sum_{k=0}^{\infty} \sum_{j=0}^{\infty} \sum_{m=0}^{\infty} \frac{(-1)^{k+j} \beta^{\alpha+k} i^m t^m}{k!\,m!\,\Gamma(\alpha)} \binom{k+\alpha-1}{j} \tau_{m,0,j-\alpha-k-1}.$$

*2.1.4 Expansion of the central moments of order m for the (1-G1)/G1 Gamma Class.*

We'll look at the development of the expansion calculations for central moments of order $m$ to the $\frac{1-G_1}{G_1}$ Gamma Class. As

$$\mu'_m = E[(X-\mu)^m] = \int_{-\infty}^{+\infty} (x-\mu)^m dF(x),$$

we have

$$\mu'_m = \sum_{r=0}^{m} \binom{m}{r}(-1)^r \mu^r \mu_{m-r}.$$

Since

$$\mu_{m-r} = \sum_{k=0}^{\infty}\sum_{j=0}^{\infty} \frac{(-1)^{k+j}\beta^{\alpha+k}}{k!\,\Gamma(\alpha)} \binom{k+\alpha-1}{j} \tau_{m-r,0,j-\alpha-k-1},$$

it follows that

$$\mu'_m = \sum_{r=0}^{m}\sum_{k=0}^{\infty}\sum_{j=0}^{\infty} \frac{(-1)^{k+j+r}\beta^{\alpha+k}\mu^r}{k!\,\Gamma(\alpha)} \binom{m}{r}\binom{k+\alpha-1}{j} \tau_{m-r,0,j-\alpha-k-1}.$$

∎

In particular, we need to expand the range of variance for the $\frac{1-G_1}{G_1}$ Gamma Class is given by:

$$\sigma^2 = \mu'_2 = \sum_{r=0}^{2}\sum_{k=0}^{\infty}\sum_{j=0}^{\infty} \frac{(-1)^{k+j+r}\beta^{\alpha+k}\mu^r}{k!\,\Gamma(\alpha)} \binom{2}{r}\binom{k+\alpha-1}{j} \tau_{2-r,0,j-\alpha-k-1}.$$

∎

*2.1.5 Expansion to the general rate for the (1-$G_1$)/$G_1$ Gamma Class*

We'll look at the development of the expansion calculations for the general coefficient for the $\frac{1-G_1}{G_1}$ Gamma Class. As

$$C_g(m) = \frac{E[(X-\mu)^m]}{\sqrt{\{E[(X-\mu)^2]\}^m}} = \frac{E[(X-\mu)^m]}{\sigma^m} = \frac{\mu'_m}{\sigma^m},$$

we have

$$C_g(m) = \frac{\sum_{r=0}^{m}\sum_{k=0}^{\infty}\sum_{j=0}^{\infty} \frac{(-1)^{k+j+r}\beta^{\alpha+k}\mu^r}{k!\,\Gamma(\alpha)} \binom{m}{r}\binom{k+\alpha-1}{j} \tau_{m-r,0,j-\alpha-k-1}}{\left(\sum_{r=0}^{2}\sum_{k=0}^{\infty}\sum_{j=0}^{\infty} \frac{(-1)^{k+j+r}\beta^{\alpha+k}\mu^r}{k!\,\Gamma(\alpha)} \binom{2}{r}\binom{k+\alpha-1}{j} \tau_{2-r,0,j-\alpha-k-1}\right)^{\frac{m}{2}}}.$$

In particular, as $C_a = C_g(3)$ we will have the expansion for the asymmetry coefficient for the $\frac{1-G_1}{G_1}$ Gamma Class is given by:

$$C_a = \frac{\sum_{r=0}^{3}\sum_{k=0}^{\infty}\sum_{j=0}^{\infty}\frac{(-1)^{k+j+r}\beta^{\alpha+k}\mu^r}{k!\Gamma(\alpha)}\binom{3}{r}\binom{k+\alpha-1}{j}\tau_{3-r,0,j-\alpha-k-1}}{\left(\sum_{r=0}^{2}\sum_{k=0}^{\infty}\sum_{j=0}^{\infty}\frac{(-1)^{k+j+r}\beta^{\alpha+k}\mu^r}{k!\Gamma(\alpha)}\binom{2}{r}\binom{k+\alpha-1}{j}\tau_{2-r,0,j-\alpha-k-1}\right)^{\frac{3}{2}}}.$$

Similarly, as $C_c = C_g(4)$, we have the expansion for the kurtosis coefficient for the $\frac{1-G_1}{G_1}$ Gamma Class is given by:

$$C_c = \frac{\sum_{r=0}^{4}\sum_{k=0}^{\infty}\sum_{j=0}^{\infty}\frac{(-1)^{k+j+r}\beta^{\alpha+k}\mu^r}{k!\Gamma(\alpha)}\binom{4}{r}\binom{k+\alpha-1}{j}\tau_{4-r,0,j-\alpha-k-1}}{\left(\sum_{r=0}^{2}\sum_{k=0}^{\infty}\sum_{j=0}^{\infty}\frac{(-1)^{k+j+r}\beta^{\alpha+k}\mu^r}{k!\Gamma(\alpha)}\binom{2}{r}\binom{k+\alpha-1}{j}\tau_{2-r,0,j-\alpha-k-1}\right)^{2}}.$$

*2. 1.6 Derivative of the log-likelihood function with respect to the parameters for the $(1-G_1)/G_1$ Gamma Class*

Once met some regularity conditions, the maximum likelihood estimators can be obtained by equating the derivative of the log-likelihood function with respect to each parameter to zero. We'll look at the calculations of the derivative of the log-likelihood function with respect to the parameters for the $\frac{1-G_1}{G_1}$ Gamma Class. As

$$\sum_{i=1}^{n} \log h_{G_1}(x_i;\alpha,\beta,\underline{\theta}) = n\text{Log}\left(\frac{\lambda\beta^\alpha}{\Gamma(\alpha)}\right) + \sum_{i=1}^{n}\log\left(\frac{g_1(x_i;\underline{\theta})}{G_1^2(x_i;\underline{\theta})}\right) +$$

$$+ (\alpha-1)\sum_{i=1}^{n}\log\left(\frac{1-G_1(x_i;\underline{\theta})}{G_1(x_i;\underline{\theta})}\right) - \beta\sum_{i=1}^{n}\left(\frac{1-G_1(x_i;\underline{\theta})}{G_1(x_i;\underline{\theta})}\right),$$

we have

$$\sum_{i=1}^{n}\frac{\partial \log h_{G_1}(x_i;\alpha,\beta,\underline{\theta})}{\partial \alpha} = n\text{Log}\beta - \frac{n\Gamma'(\alpha)}{\Gamma(\alpha)} + \sum_{i=1}^{n}\log\left(\frac{1-G_1(x_i;\underline{\theta})}{G_1(x_i;\underline{\theta})}\right).$$

$$\sum_{i=1}^{n}\frac{\partial \log h_{G_1}(x_i;\alpha,\beta,\underline{\theta})}{\partial \beta} = \frac{n\alpha}{\beta} - \sum_{i=1}^{n}\left(\frac{1-G_1(x_i;\underline{\theta})}{G_1(x_i;\underline{\theta})}\right).$$

$$\sum_{i=1}^{n}\frac{\partial \log h_{G_1}(x_i;\alpha,\beta,\underline{\theta})}{\partial \theta_j} = \sum_{i=1}^{n}\frac{\partial \log\left(\frac{g_1(x_i;\underline{\theta})}{G_1^2(x_i;\underline{\theta})}\right)}{\partial \theta_j} + (\alpha-1)\sum_{i=1}^{n}\frac{\partial \log\left(\frac{1-G_1(x_i;\underline{\theta})}{G_1(x_i;\underline{\theta})}\right)}{\partial \theta_j} - \beta\sum_{i=1}^{n}\frac{\partial \left(\frac{1-G_1(x_i;\underline{\theta})}{G_1(x_i;\underline{\theta})}\right)}{\partial \theta_j}.$$

∎

*2.1.7 Entropy Rényi using the $(1-G_1)/G_1$ Gamma Class*

Entropy is a measure of uncertainty in the sense that the higher the entropy value the lowest the information and the greater the uncertainty, or the greater the randomness or disorder. The following is the expansion entropy calculations for the $\frac{1-G_1}{G_1}$ Gamma Class, using the Rényi *entropy*, which is given by

$$L_R(\eta) = \frac{1}{1-\eta} \log \left( \int_{-\infty}^{+\infty} f^\eta(x) dF(x) \right).$$

Substituting the expressions of density and cumulative distribution function, we have

$$L_R(\eta) = \frac{1}{1-\eta} \log \left( \int_{-\infty}^{+\infty} \left( \frac{g_1(x)}{G_1^2(x)} \frac{\beta^\alpha}{\Gamma(\alpha)} \left( \frac{1-G_1(x)}{G_1(x)} \right)^{\alpha-1} \exp\left( -\beta \left( \frac{1-G_1(x)}{G_1(x)} \right) \right) \right)^\eta dx \right).$$

As

$$\exp\left( -\eta\beta \frac{1-G_1(x)}{G_1(x)} \right) = \sum_{k=0}^{\infty} \frac{(-1)^k \eta^k \beta^k}{k!} G_1^{-k} (1-G_1(x))^k,$$

we have

$$h_{G_1}^\eta(x) = g_1^\eta(x) \frac{\beta^{\eta\alpha}}{\Gamma^\eta(\alpha)} \sum_{k=0}^{\infty} \frac{(-1)^k \eta^k \beta^k}{k!} G_1^{-\eta(\alpha+1)-k}(x) (1-G_1(x))^{\eta(\alpha-1)+k}.$$

Using the following expansion

$$(1-G_1(x))^{\eta(\alpha-1)+k} = \sum_{j=0}^{\infty} \binom{\eta(\alpha-1)+k}{j} (-1)^j G_1^j(x),$$

it follows that

$$h_{G_1}^\eta(x) = \sum_{k=0}^{\infty} \sum_{j=0}^{\infty} \frac{(-1)^{k+j} \eta^k \beta^{\eta\alpha+k}}{k! \Gamma^\eta(\alpha)} \binom{\eta(\alpha-1)+k}{j} g_1^\eta(x) G_1^{-\eta(\alpha+1)-k+j}.$$

Thus, we have

$$L_R(\eta) = \frac{1}{1-\eta} \log \left( \sum_{k=0}^{\infty} \sum_{j=0}^{\infty} \frac{(-1)^{k+j} \eta^k \beta^{\eta\alpha+k}}{k! \Gamma^\eta(\alpha)} \binom{\eta(\alpha-1)+k}{j} \int_{-\infty}^{+\infty} g_1^\eta(x) G_1^{-\eta(\alpha+1)-k+j} dx \right),$$

which, in turn, implies that

$$L_R(\eta) = \frac{1}{1-\eta} \log \left( \sum_{k=0}^{\infty} \sum_{j=0}^{\infty} \frac{(-1)^{k+j} \eta^k \beta^{\eta\alpha+k}}{k! \Gamma^\eta(\alpha)} \binom{\eta(\alpha-1)+k}{j} \tau_{0,\eta-1,-\eta(\alpha+1)-k+j} \right).$$

## 3. Distribution proposal

### 3.1 Construction of a distribution from the (1-G1)/G1 Gamma Class

In this section, we will examine a particular distribution of the $\frac{1-G_1}{G_1}$Gamma Class proposed here. We will consider the particular case in which $G_1(x) = 1 - e^{-\lambda x}, x > 0$, that is called the *(1-Exp)/Exp* Gamma Distribution.

#### 3.1.1 (1-Exp) / Exp Gamma distribution

Considering $G_1(x)$ the cdf of the exponential distribution with parameter $\lambda$ in the generator functional for the $\frac{1-G_1}{G_1}$Gamma Class, we have the *(1-Exp)/Exp* Gamma distribution:

$$H(x) = 1 - \int_0^{\frac{e^{-\lambda x}}{1-e^{-\lambda x}}} \frac{\beta^\alpha}{\Gamma(\alpha)} t^{\alpha-1} e^{-\beta t} dt, x > 0,$$

which is equivalently given by

$$H(x) = \int_{\frac{e^{-\lambda x}}{1-e^{-\lambda x}}}^{+\infty} \frac{\beta^\alpha}{\Gamma(\alpha)} t^{\alpha-1} e^{-\beta t} dt, x > 0.$$

Differentiating $H(x)$, we get the density function of the $\frac{1-exp}{exp}$ Gamma distribution:

$$h(x) = \frac{\lambda \beta^\alpha}{\Gamma(\alpha)} \frac{e^{-\lambda x}}{(1-e^{-\lambda x})^2} \left(\frac{e^{-\lambda x}}{1-e^{-\lambda x}}\right)^{\alpha-1} e^{-\beta\left(\frac{e^{-\lambda x}}{1-e^{-\lambda x}}\right)}. \blacksquare$$

Figures 3.1.1 to 3.1.6 show the graphs of the $\frac{1-exp}{exp}$ Gamma distribution probability density functions and cumulative distribution, for some values of the parameters, visualizing the variation of one as a function of the other two fixed.

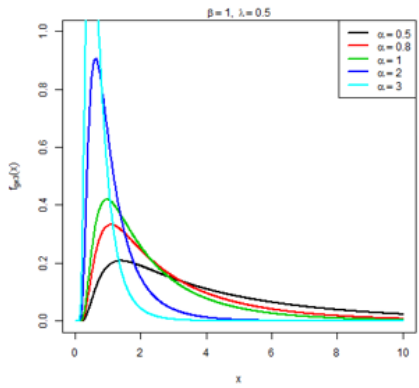

Figure 3.1.1 – $pdf$ of the gamma distribution $\frac{1-Exp}{Exp}$ with $\alpha$ variand

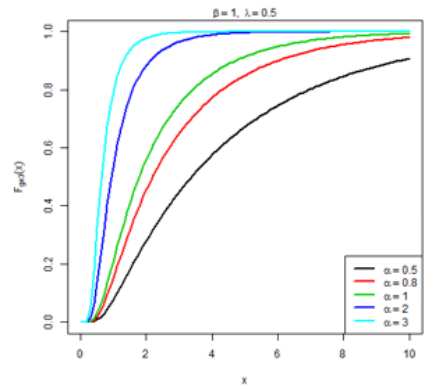

Figura 3.1.4 – $cdf$ of the gamma distribution $\frac{1-Exp}{Exp}$ with $\alpha$ variand

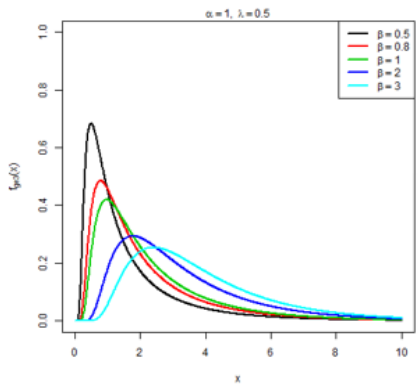

Figure 3.1.2 – $pdf$ of the gamma distribution $\frac{1-Exp}{Exp}$ with $\beta$ variand

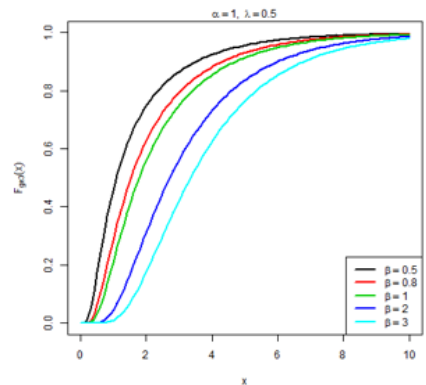

Figure 3.1.5 – $cdf$ of the gamma distribution $\frac{1-Exp}{Exp}$ com $\beta$ variand

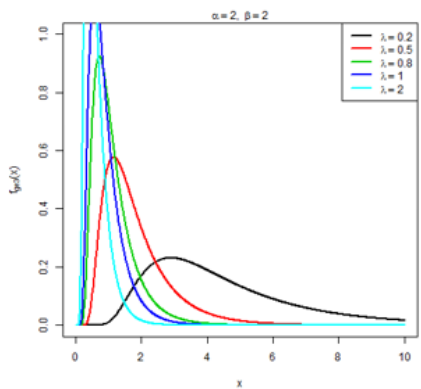

Figure 3.1.3 – $pdf$ of the gamma distribution $\frac{1-Exp}{Exp}$ com $\lambda$ variand

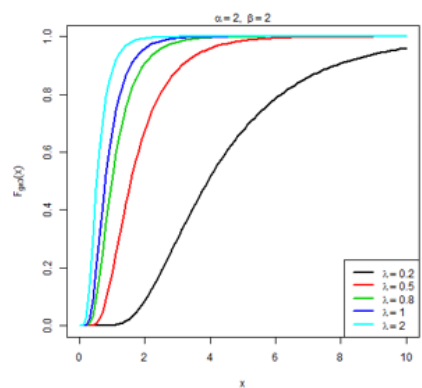

Figure 3.1.6 $cdf$ of the gamma distribution $\frac{1-Exp}{Exp}$ com $\lambda$ variand

Using procedure similar to what was done in Section 2, we can rewrite the density and cumulative distribution function of the gamma distribution as a sum of exponentiated Exponentials, as follows:

$$h(x) = \sum_{k=0}^{\infty}\sum_{j=0}^{\infty} \binom{-k-\alpha-1}{j} \frac{(-1)^{k+j}\lambda\beta^{k+\alpha}}{k!\,\Gamma(\alpha)} e^{-\lambda(k+\alpha+j)x}$$

and

$$H(x) = \sum_{k=0}^{\infty}\sum_{j=0}^{\infty} \binom{-k-\alpha-1}{j} \frac{(-1)^{k+j+1}\beta^{k+\alpha}}{k!\,(k+\alpha+j)\Gamma(\alpha)} \left(e^{-\lambda(k+\alpha+j)x}-1\right).$$

*3.1.2 Risk function using the (1-Exp) / Exp Gamma distribution*

We can also obtain the risk function using the $\frac{1-exp}{exp}$ Gamma distribution as follows:

$$R(x) = \frac{h(x)}{1-H(x)}.$$

$$R(x) = \frac{\dfrac{\lambda\beta^{\alpha}}{\Gamma(\alpha)}\dfrac{e^{-\lambda x}}{(1-e^{-\lambda x})^2}\left(\dfrac{e^{-\lambda x}}{1-e^{-\lambda x}}\right)^{\alpha-1} e^{-\beta\left(\frac{e^{-\lambda x}}{1-e^{-\lambda x}}\right)}}{1 - \int_{\frac{e^{-\lambda x}}{1-e^{-\lambda x}}}^{+\infty} \dfrac{\beta^{\alpha}}{\Gamma(\alpha)} t^{\alpha-1}e^{-\beta t}dt}$$

Therefore,

$$R(x) = \frac{\dfrac{\lambda\beta^{\alpha}}{\Gamma(\alpha)}\dfrac{e^{-\lambda x}}{(1-e^{-\lambda x})^2}\left(\dfrac{e^{-\lambda x}}{1-e^{-\lambda x}}\right)^{\alpha-1} e^{-\beta\left(\frac{e^{-\lambda x}}{1-e^{-\lambda x}}\right)}}{\int_0^{\frac{e^{-\lambda x}}{1-e^{-\lambda x}}} \dfrac{\beta^{\alpha}}{\Gamma(\alpha)} t^{\alpha-1}e^{-\beta t}dt}$$

Figures 3.2.1 to 3.2.3 show the graphs of the risk function using the $\frac{1-exp}{exp}$ Gamma distribution generated from some values assigned to parameters.

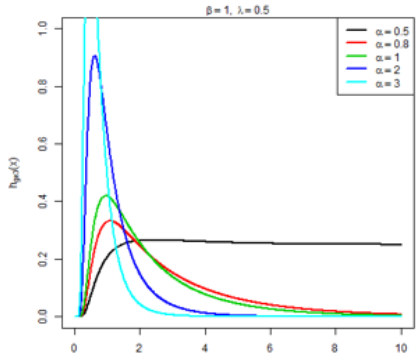
Figure 3.2.1 - $\mathcal{R}(x)$ of the distribution range $\frac{1-Exp}{Exp}$ with varying $\alpha$

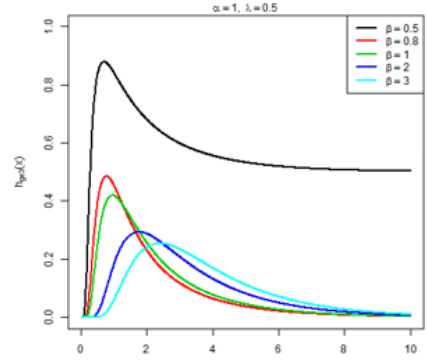
Figure 3.2.2 - $\mathcal{R}(x)$ of the distribution range $\frac{1-Exp}{Exp}$ with varying $\beta$

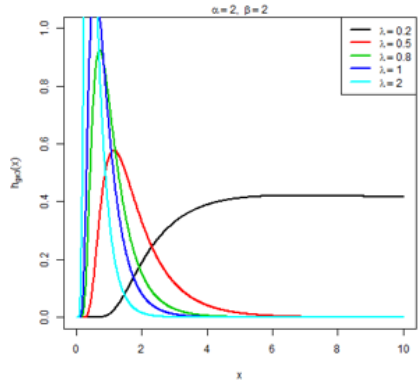
Figure 3.2.3 - $\mathcal{R}(x)$ of the distribution range $\frac{1-Exp}{Exp}$ with varying $\lambda$

*3.1.3 Other properties of the Gamma distribution (1-Exp) / Exp*

Using the $pdf$ and $cdf$ expansion, we can obtain the statistical properties of the $\frac{1-exp}{exp}$ Gamma distribution using procedures similar to those used in Section 2 to calculate the corresponding properties of $\frac{1-G_1}{G_1}$ Gamma Class. In this way, we obtain the following expressions:

Expansion for times of order m

$$\mu_m = \sum_{k=0}^{\infty}\sum_{j=0}^{\infty} \binom{-k-\alpha-1}{j} \frac{(-1)^{k+j}\lambda\beta^{k+\alpha}}{k!\,\Gamma(\alpha)} \frac{\Gamma(m+1)}{\left(\lambda(k+\alpha+j)\right)^{m+1}}$$

In particular we have that the mean of the $\frac{1-exp}{exp}$Gamma distribution:

$$\mu_1 = \mu = \sum_{k=0}^{\infty}\sum_{j=0}^{\infty}\binom{-k-\alpha-1}{j}\frac{(-1)^{k+j}\lambda\beta^{k+\alpha}}{k!\Gamma(\alpha)}\frac{1}{(\lambda(k+\alpha+j))^2}.$$

*Expansion for the moment generating function*

$$M_X(t) = \sum_{k=0}^{\infty}\sum_{j=0}^{\infty}\binom{-k-\alpha-1}{j}\frac{(-1)^{k+j}\lambda\beta^{k+\alpha}}{k!\Gamma(\alpha)}\frac{1}{\lambda(k+\alpha+j)-t}, \text{ if } t < \lambda(k+\alpha+j)$$

Expansion for the characteristic function

$$\varphi_X(t) = \sum_{k=0}^{\infty}\sum_{j=0}^{\infty}\binom{-k-\alpha-1}{j}\frac{(-1)^{k+j}\lambda\beta^{k+\alpha}}{k!\Gamma(\alpha)}\frac{1}{\lambda(k+\alpha+j)-it}.$$

Expansion for central moments of order m

$$\mu'_m = \sum_{r=0}^{m}\sum_{k=0}^{\infty}\sum_{j=0}^{\infty}\binom{-k-\alpha-1}{j}\binom{m}{r}\frac{(-1)^{k+j+r}\lambda\beta^{k+\alpha}}{k!\Gamma(\alpha)}\frac{\mu^r\Gamma(m-r+1)}{(\lambda(k+\alpha+j))^{m-r+1}}.$$

In particular, the expansion to the variance of the $\frac{1-exp}{exp}$Gamma distribution *is given* by:

$$\sigma^2 = \mu'_2 = \sum_{r=0}^{2}\sum_{k=0}^{\infty}\sum_{j=0}^{\infty}\binom{-k-\alpha-1}{j}\binom{2}{r}\frac{(-1)^{k+j+r}\lambda\beta^{k+\alpha}}{k!\Gamma(\alpha)}\frac{\mu^r\Gamma(3-r)}{(\lambda(k+\alpha+j))^{3-r}}.$$

Expansion to the general coefficient

$$C_g(m) = \frac{\sum_{r=0}^{m}\sum_{k=0}^{\infty}\sum_{j=0}^{\infty}\binom{-k-\alpha-1}{j}\binom{m}{r}\frac{(-1)^{k+j+r}\lambda\beta^{k+\alpha}}{k!\Gamma(\alpha)}\frac{\mu^r\Gamma(m-r+1)}{(\lambda(k+\alpha+j))^{m-r+1}}}{\left(\sum_{r=0}^{2}\sum_{k=0}^{\infty}\sum_{j=0}^{\infty}\binom{-k-\alpha-1}{j}\binom{2}{r}\frac{(-1)^{k+j+r}\lambda\beta^{k+\alpha}}{k!\Gamma(\alpha)}\frac{\mu^r\Gamma(3-r)}{(\lambda(k+\alpha+j))^{3-r}}\right)^{\frac{m}{2}}}.$$

In particular, the expansion for the coefficient of asymmetry for the $\frac{1-exp}{exp}$Gamma distribution, $C_a = C_g(3)$, is given by:

$$C_a = \frac{\sum_{r=0}^{3} \sum_{k=0}^{\infty} \sum_{j=0}^{\infty} \binom{-k-\alpha-1}{j}\binom{3}{r}\frac{(-1)^{k+j+r}\lambda\beta^{k+\alpha}}{k!\,\Gamma(\alpha)}\frac{\mu^r\Gamma(4-r)}{(\lambda(k+\alpha+j))^{4-r}}}{\left(\sum_{r=0}^{2} \sum_{k=0}^{\infty} \sum_{j=0}^{\infty} \binom{-k-\alpha-1}{j}\binom{2}{r}\frac{(-1)^{k+j+r}\lambda\beta^{k+\alpha}}{k!\,\Gamma(\alpha)}\frac{\mu^r\Gamma(3-r)}{(\lambda(k+\alpha+j))^{3-r}}\right)^{\frac{3}{2}}}.$$

While the kurtosis expansion coefficient of the $\frac{1-\exp}{\exp}$ Gamma distribution, $C_c = C_g(4)$, is given by:

$$C_c = \frac{\sum_{r=0}^{4} \sum_{k=0}^{\infty} \sum_{j=0}^{\infty} \binom{-k-\alpha-1}{j}\binom{4}{r}\frac{(-1)^{k+j+r}\lambda\beta^{k+\alpha}}{k!\,\Gamma(\alpha)}\frac{\mu^r\Gamma(5-r)}{(\lambda(k+\alpha+j))^{5-r}}}{\left(\sum_{r=0}^{2} \sum_{k=0}^{\infty} \sum_{j=0}^{\infty} \binom{-k-\alpha-1}{j}\binom{2}{r}\frac{(-1)^{k+j+r}\lambda\beta^{k+\alpha}}{k!\,\Gamma(\alpha)}\frac{\mu^r\Gamma(3-r)}{(\lambda(k+\alpha+j))^{3-r}}\right)^{2}}.$$

Derivatives from the log-likelihood function with respect to the parameters

As

$$\sum_{i=1}^{n} Log\,h(x_i;\alpha,\beta,\lambda) = nLog\left(\frac{\lambda\beta^\alpha}{\Gamma(\alpha)}\right) - n\alpha\lambda - (\alpha+1)\sum_{i=1}^{n} \log(1-e^{-\lambda x_i}) - \beta\sum_{i=1}^{n}\left(\frac{e^{-\lambda x_i}}{1-e^{-\lambda x_i}}\right),$$

we have

$$\sum_{i=1}^{n} \frac{\partial Log\,h(x_i;\alpha,\beta,\lambda)}{\partial \alpha} = nLog\beta - \frac{n\Gamma'(\alpha)}{\Gamma(\alpha)} - n\lambda - \sum_{i=1}^{n} \log(1-e^{-\lambda x_i}).$$

$$\sum_{i=1}^{n} \frac{\partial Log\,h(x_i;\alpha,\beta,\lambda)}{\partial \beta} = \frac{n\alpha}{\beta} - \sum_{i=1}^{n}\left(\frac{e^{-\lambda x_i}}{1-e^{-\lambda x_i}}\right).$$

$$\sum_{i=1}^{n} \frac{\partial Log\,h(x_i;\alpha,\beta,\lambda)}{\partial \lambda} = \frac{n}{\lambda} - n\alpha - (\alpha+1)\sum_{i=1}^{n}\left(\frac{x_i e^{-\lambda x_i}}{1-e^{-\lambda x_i}}\right) + \beta\sum_{i=1}^{n}\left(\frac{x_i e^{-\lambda x_i}}{(1-e^{-\lambda x_i})^2}\right).$$

Rényi entropy

$$L_R(\eta) = \frac{1}{1-\eta}\log\left(\sum_{k=0}^{\infty}\sum_{j=0}^{\infty}\frac{(-1)^{k+j}\lambda^\eta\eta^k\beta^{\eta\alpha+k}}{k!\,\Gamma^\eta(\alpha)}\binom{\eta(\alpha-1)+k}{j}\frac{1}{\lambda(\alpha+k+j)}\right).$$

*3.1.4 Application*

In this section, we will show an application to real data for the proposed Gamma distribution. The data used in this research are from the excesses of flood peaks (in m $^3$ / s) Wheaton River near Carcross in the Yukon Territory, Canada. 72 exceedances of the years 1958 to 1984 were recorded, rounded to one decimal place. These data were analyzed by Choulakian and Stephens (2001), and are presented in Table 3.1.1.

Table 3.1.1: full excess peaks in m $^3$ / s Rio Wheaton

| Excess flood peaks of Rio Wheaton (m $^3$ / s) | | | | | | | | | | | |
|---|---|---|---|---|---|---|---|---|---|---|---|
| 1.7 | 2.2 | 14.4 | 1.1 | 0.4 | 20.6 | 5.3 | 0.7 | 1.9 | 13.0 | 12.0 | 9.3 |
| 1.4 | 18.7 | 8.5 | 25.5 | 11.6 | 14.1 | 22.1 | 1.1 | 2.5 | 14.4 | 1.7 | 37.6 |
| 0.6 | 2.2 | 39.0 | 0.3 | 15.0 | 11.0 | 7.3 | 22.9 | 1.7 | 0.1 | 1.1 | 0.6 |
| 9.0 | 1.7 | 7.0 | 20.1 | 0.4 | 2.8 | 14.1 | 9.9 | 10.4 | 10.7 | 30.0 | 3.6 |
| 5.6 | 30.8 | 13.3 | 4.2 | 25.5 | 3.4 | 11.9 | 21.5 | 27.6 | 36.4 | 2.7 | 64.0 |
| 1.5 | 2.5 | 27.4 | 1.0 | 27.1 | 20.2 | 16.8 | 5.3 | 9.7 | 27.5 | 2.5 | 27.0 |

It is worth mentioning that this data set has also been analyzed by means of the distributions of Pareto, Weibull three parameters, the generalized Pareto and beta - Pareto (AKINSETE *et al*, 2008).

In Table 3.1.2, we can see the maximum likelihood estimates obtained by the Newton-Raphson implemented in SAS 9.1 statistical software, parameters, standard errors, Akaike information criterion, corrected Akaike, Bayesian, Hannan-Quinn and Anderson-Darling statistics (A) and Cramér von Mises (W) to the $-ln(1-exp)$ Gamma distributions *(M 1)*, $\frac{1-exp}{exp}$ Gamma distribution *(proposed model, M 2)*, exponentiated

Weibull (M3), modified Weibull (M4), beta Pareto (M5) and Weibull (M6).

Table 3.1.2 - Estimated maximum likelihood parameter, errors (standard errors in parentheses) and calculations of AIC statistics, AIC, BIC, HQIC, tests A and W for the M1 to M6 distributions.

| Models | Statistics | | | | Statistics | | | | | |
|---|---|---|---|---|---|---|---|---|---|---|
| | | | | | AIC | AIC | BIC | HQIC | The | W |
| M1 | 0.838 (0.121) | 0.035 (0.007) | 1.96 (<E-3) | ---- ---- | 508.689 | 509.042 | 515.519 | 511.408 | 0.7519 | 0.1306 |
| M2 | 0.131 (0.053) | 0.179 (0.07) | 0.539 (0.251) | ---- ---- | 505.030 | 505.383 | 511.860 | 507.749 | 0.4516 | 0.0757 |
| M3 | 1,387 (0.59) | 0.519 (0.312) | 0.05 (0.021) | ---- ---- | 508.050 | 508.403 | 514.880 | 510.769 | 1.4137 | 0.2534 |
| M4 | 0.776 (0.124) | 0.124 (0.035) | 0.01 (0.008) | ---- ---- | 507.343 | 507.696 | 514.173 | 510.062 | 0.5938 | 0.0978 |
| M5 | 84.682 (<E-3) | 65.574 (<E-3) | 0.063 (0.005) | 0.01 (<E-3) | 524.398 | 524.995 | 533.504 | 528.023 | 2.0412 | 0.3516 |
| M6 | 0.901 (0.086) | 0.086 (0.012) | ---- ---- | ---- ---- | 506.997 | 507.171 | 511.551 | 508.810 | 0.7855 | 0.1380 |

For the six distributions shown in Table 3.1.2, the data applied to Wheaton River flooding, it was observed that beta-Pareto model (M5), which was described in

2008 by Akinsete, Famoye and Lee as the best fitted model, in our studies had a lower performance with AIC = 524.398, AICc = 524.995, BIC=533.504, HQIC = 528.023, A = 2.0412 and W = 0.3516, when compared to the proposed $\frac{1-exp}{exp}$Gamma model *(M 2)* that obtained AIC = 505.030, AICc = 505.383, BIC=511.860, HQIC = 507.749, A = 0.4516 and W = 0.0757. Also according to Table 3.1.2, the proposed distribution model M2 is the best tested once the loweest values of AIC, AICc, BIC HQIC, A and W are from such distribution, and only according to the BIC criterion was exceeded solely by the model M6.

In the Figures 3.3.1 and 3.3.2 below, there are the graphs of density functions and distributions of M1 to M6 models fitted to the data and their corresponding histograms. The graph shows that the $\frac{1-exp}{exp}$Gamma model has similar behavior to that of other distributions, except that of the beta Pareto which distances itself from the others.

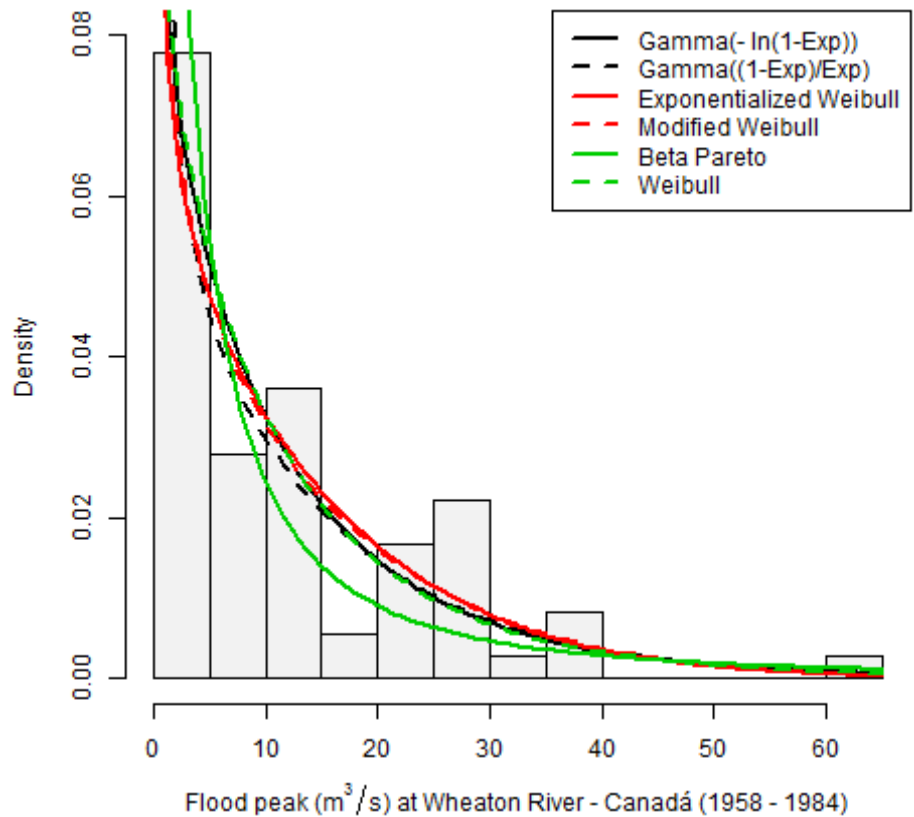

Figure 3.3.1 adjusted to the mass data of flood peaks in river Wheaton

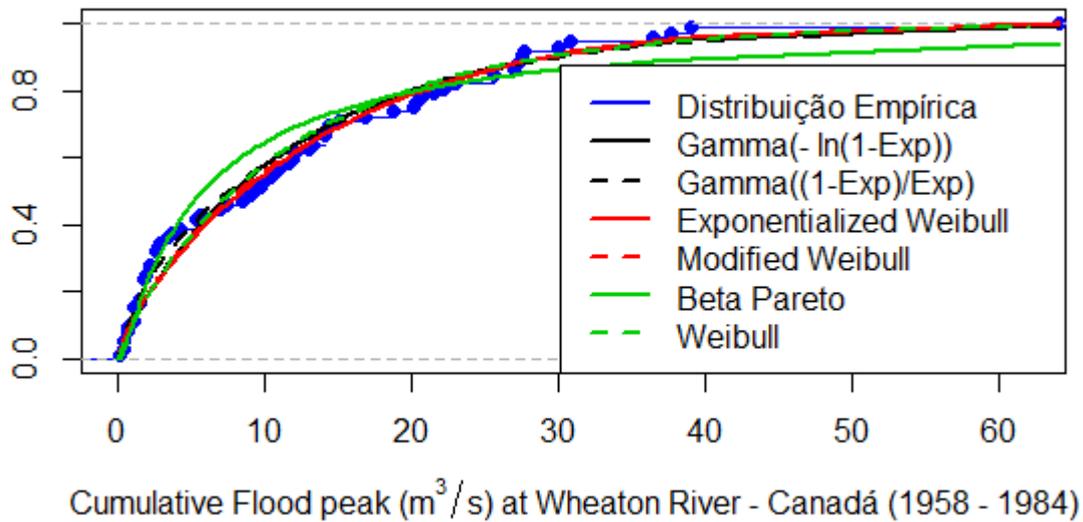

Figure 3.3.2 adjusted to the mass data of flood peaks in river Wheaton

## 4. Conclusion

As concluding remarks, we note that the class of $(1-G_1)/G_1$ gamma probability distributions developed in this work is a novel way of generalizing the gamma distribution and can be applied in different areas depending on the choice of the distribution $G_1$. As future work, we intend to carry out more detailed comparisons between the novel distribution family proposed in this paper and the family of distributions investigated in Zografos and Balakrishnan (2009) and Silva (2013), which are also based on the integration of the gamma distribution.

In this work, we study in detail only a distribution of the $\frac{1-G_1}{G_1}$ Gamma Class, namely the $\frac{1-exp}{exp}$ Gamma distribution. We derive the properties of this distribution and applied to a set of real data obtaining better fit than that obtained in a previous study by Akinsete *et al.* (2008). We intend to conduct the study of new distributions within this class as future work.

We note that after adding several parameters to a model it can better be adjusted to a particular phenomenon due to its greater flexibility. On the other hand, one should not forget that there may be a problem for the estimation of the parameters since it can occur both computational and identifiability problems in parameter estimation. Thus, the ideal is to choose a model that reflects well the phenomenon / experiment with the minimum number of parameters. In the case of the proposed class, only two additional parameters are added to the set of parameters of the G distribution.